\documentclass[11pt,reqno]{amsart}
\usepackage{amssymb}
\usepackage{palatino}
\input amssym.def
\usepackage{amsmath, amsfonts}
\usepackage{amssymb}
\usepackage{amscd, mathtools}
\usepackage[mathscr]{eucal}
\usepackage{palatino}
\newfont{\cyrr}{wncyr10}

\newcommand{\Z}{{\mathbb Z}}
\newcommand{\Q}{{\mathbb Q}}

\newcommand{\C}{{\mathbb C}}

\newcommand{\K}{{\mathbf K}}
\newcommand{\E}{{\mathbf E}}
\newcommand{\F}{{\mathbf F}}
\renewcommand{\L}{{\mathbf L}}

\renewcommand{\mod}{{\, \rm mod \, }}
\newcommand{\thmref}[1]{Theorem~\ref{#1}}
\newtheorem{thm}{Theorem}

\newtheorem{lem}[thm]{Lemma}
\newtheorem{cor}[thm]{Corollary}
\newtheorem{prop}[thm]{Proposition}
\newtheorem{rmk}{Remark}[section] 
\newtheorem{defn}{Definition}
\newtheorem{conj}{Conjecture}

\newcommand{\corref}[1]{Corollary~\ref{#1}}
\newcommand{\propref}[1]{Proposition~\ref{#1}}
\newcommand{\lemref}[1]{Lemma~\ref{#1}}

\begin{document}

\title[Extension of a question of Baker]{On an extension of a question of Baker}

\author{sanoli Gun and Neelam Kandhil }

\address{Sanoli Gun and Neelam Kandhil \\ \newline
The Institute of Mathematical Sciences, A CI of Homi Bhabha National Institute, 
CIT Campus, Taramani, Chennai 600 113, India.}
\email{sanoli@imsc.res.in} 
\email{neelam@imsc.res.in}



\maketitle   

\begin{abstract} 
It is an open question of Baker whether the numbers $L(1, \chi)$
for non-trivial Dirichlet characters $\chi$ with period $q$ are linearly independent 
over $\Q$. The best known result is due to Baker, Birch and Wirsing
which affirms this when $q$ is co-prime to $\varphi(q)$.  In this article, 
we extend their result to any arbitrary family of moduli.
More precisely, for a positive integer $q$, let $X_q$
denote the set of all $L(1,\chi)$ values
as $\chi$ varies over non-trivial Dirichlet characters with period $q$.
 Then for any finite set of pairwise co-prime natural numbers $q_i, 1\le i \le \ell$ 
with $(q_1 \cdots q_{\ell}, ~\varphi(q_1)\cdots \varphi(q_{\ell}))=1$,
we show that the set $X_{q_1} \cup \cdots \cup X_{q_l}$
is linearly independent over $\Q$.
In the process, we also extend a result of Okada about linear independence 
of the cotangent values over $\Q$ as well as a result of Murty-Murty about
$\overline{\Q}$ linear independence of such $L(1, \chi)$ values. Finally, 
we prove $\Q$ linear independence of such $L$ values of 
Erd\"{o}sian  functions with distinct prime periods $p_i$ for $1\le i \le \ell$ with 
$(p_1 \cdots p_{\ell}, ~ \varphi( p_1\cdots p_{\ell}) )=1$. 
\end{abstract}

\section{\bf{Introduction}}
For an integer $q>1$ and a Dirichlet character $\chi$ with period $q$, consider the
Dirichlet $L$ function
$$
L(s, \chi) = \sum_{n=1}^{\infty} \frac{\chi(n)}{n^s}, \phantom{m} \Re(s) > 1.
$$
When $\chi$ is non-trivial, we know that $L(s,\chi)$
extends to an entire function, $L(1, \chi)$ is non-zero and is
equal to $\sum_{n=1}^{\infty} \frac{\chi(n)}{n}$.
For $q$ as before, consider the set
$$
X_q = \{ L(1,\chi) ~|~ \chi \mod{q},  \chi \ne \chi_0 \},
$$
where $\chi_0$ is the trivial Dirichlet character with period $q$.
In \cite[p. 48]{AB}, Baker asked whether the numbers
in $X_q$ are linearly independent over $\Q$.
In 1973, Baker, Birch and Wirsing \cite{BBW} in an elegant work
proved that the numbers in the set 
$\{ L(1,\chi) \in X_q ~|~   (q, \varphi(q))=1 \}$
are linearly independent over $\Q$ (see \cite{NR} for
an exposition on this topic).
In this context, we prove the following theorem.

\begin{thm}\label{th1}
Let $q_j >2$ for $1 \le j \le \ell$ be pairwise co-prime natural numbers
such that the number $q_1\cdots q_{\ell}$ is co-prime to $\varphi(q_1)\cdots \varphi(q_{\ell}))$. 
Then the numbers in the set 
$X_{q_1} \cup \cdots \cup X_{q_\ell} $
are linearly independent over $\Q(\zeta_{\varphi(q_1)\cdots \varphi(q_{\ell})})$. 
\end{thm}

More generally, we derive the following theorem.

\begin{thm}\label{th2}
Let $q_j >2$ for $1 \le j \le \ell$ be pairwise co-prime natural numbers.
Also let $\K$ be a number field with 
$\K(\zeta_{\varphi(q_1)\cdots \varphi(q_{\ell})})\cap \Q(\zeta_{q_1\cdots q_{\ell} }) = \Q$,
where $\zeta_q$ denotes a primitive $q$th root of unity.
Then the numbers in the set $X_{q_1} \cup \cdots \cup X_{q_\ell}$
are linearly independent over $\K(\zeta_{\varphi(q_1)\cdots \varphi(q_{\ell})})$. 
\end{thm}
Note that $X_q = X_{q, e} \cup X_{q, o}$, where
\begin{eqnarray}\label{eos}
X_{q, e} 
&=& \{ L(1,\chi) ~|~ \chi \mod{q}, ~ \chi(-1) = 1, \chi \ne \chi_0 \} \nonumber \\
\phantom{m}\text{and}\phantom{m}
X_{q, o} 
&=& \{ L(1,\chi) ~|~ \chi \mod{q},  ~\chi(-1) = -1 \}.\nonumber
\end{eqnarray}
In 2011, Murty-Murty refined Baker-Birch-Wirsing result to show that

\begin{thm}{\rm (Murty-Murty \cite{MM1})}
Let $q>2$ be a natural number. Then the numbers in the set 
$X_{q, e}$ are linearly independent over $\overline{\Q}$. 
\end{thm}
In this article, we prove the following theorem.

\begin{thm}\label{th3}
Let $q_j >2$ for $1 \le j \le \ell$ be pairwise co-prime natural numbers.
Then the numbers in the set $X_{q_1, e} \cup \cdots \cup X_{q_{\ell}, e}$
are linearly independent over $\overline{\Q}$. 
\end{thm}

In 1981, Okada \cite{TO} (see also \cite{ {MS1}, {KW}}) proved that
\begin{thm}{\rm(Okada \cite{TO})}
Let $q>2$ be a natural number and $\K$ be a number field with the property that
$\K(\zeta_{\varphi(q)}) \cap \Q(\zeta_q) = \Q$. Then the numbers in the set 
$X_{q, o}$
are linearly independent over $\K(\zeta_{\varphi(q)})$. 
\end{thm}
Here we prove the following theorem.

\begin{thm}\label{th4}
For $1 \le j \le \ell$, let $q_j >2$  be pairwise co-prime natural numbers.
If $\K$ is a number field such that 
$\K(\zeta_{\varphi(q_1)\cdots \varphi(q_{\ell})}) \cap \Q(\zeta_{q_1\cdots q_{\ell} }) = \Q$,
then the numbers in the set $X_{q_1,o} \cup \cdots \cup X_{q_\ell, o}$
are linearly independent over $\K(\zeta_{\varphi(q_1)\cdots \varphi(q_{\ell})})$. 
\end{thm}

\begin{rmk}
In \thmref{th2} and \thmref{th4}, $\varphi(q_1)\cdots \varphi(q_{\ell})$
in $\K(\zeta_{\varphi(q_1)\cdots \varphi(q_{\ell})})$ and 
 $q_1\cdots q_{\ell} $  in the
$ \Q(\zeta_{q_1\cdots q_{\ell} }) $
can be replaced by their respective least common multiples.
\end{rmk}

As consequences to \thmref{th2}, \thmref{th3} and \thmref{th4},
we derive the following corollaries. Before we state the corollaries,
let us introduce the notion of Dirichlet type functions as defined 
by Murty and Saradha (see \cite{ {MS1}}).
\begin{defn}
An arithmetical function $f$ with period $q>1$ with values in $\overline{\Q}$ is called 
Dirichlet type if $f(a)=0$ whenever $(a, q) \ne 1$. 
\end{defn}
\begin{defn}
A periodic function $f$ with period $q>1$ is called an Erd\"osian  function  
if $f(a)=\pm1$ for all $1 \le a <q$ and $f(q)=0$.
\end{defn}

For an arithmetical function $f$ with period $q$, consider the series
$L(s,f) = \sum_{n=1}^{\infty} \frac{f(n)}{n^s}$ for $\Re(s) >1$. 
This series has a meromorphic continuation to $\C$ with
a possible simple pole at $s=1$ of residue $q^{-1}\sum_{a=1}^{q} f(a)$
(see \cite[Ch 22]{MR} for further details). 
From now onwards, we assume that $\sum_{a=1}^{q} f(a) =0$.

For a natural number $q>2$ and a number field $\K$, let
$Y_q(\K)$ be $\K$ linearly independent set of Dirichlet type functions
of period $q$. Define
\begin{eqnarray*}
X_{q} (\K) 
&=& \{ L(1, f) ~|~  f \in Y_q(\K) \}  \\
\phantom{m}\text{and}\phantom{m}
X_{q, e}(\K) 
&=& \{ L(1, f) ~|~  f \in Y_q(\K) , ~ f(-a) = f(a) \text{ for } 1 \le a <q \}.
\end{eqnarray*}

In this set-up, we have the following corollaries. 

\begin{cor}\label{cor1}
For $1\le j \le \ell$, let $q_j >2$ be pairwise co-prime natural numbers.
Then the numbers in the set  
$X_{q_1, e}(\overline{\Q}) \cup \cdots \cup X_{q_{\ell, e}}(\overline{\Q}) $
are $\overline{\Q}$ linearly independent. 
\end{cor}

\begin{cor}\label{cor3}
 For any odd prime $p$, choose an Erd\"{o}sian function 
 $f_p$ with period $p$ which is not an
odd function. Then the numbers in the set $\{ L(1, f_p) ~|~ p \text{ odd} \}$
are linearly independent over $\overline{\Q}$.
\end{cor}

\begin{cor}\label{cor2}
Let $q_j >2, 1\le j \le \ell$ be pairwise co-prime natural numbers.
Also let $f_j  \in Y_{q_j}(\K)$ with values in a number field $\K$.
If $\K(\zeta_{\varphi(q_1)\cdots \varphi(q_{\ell})}) \cap \Q(\zeta_{q_1\cdots q_{\ell} } ) = \Q$, then
the elements in $X_{q_1}(\K) \cup \cdots \cup X_{q_{\ell}}(\K)$
are $\K$ linearly independent. In particular, choose
Erd\"{o}sian functions $f_{p_i}$ with odd prime periods $p_i$,
then the numbers $L(1, f_{p_i})$ for $1 \le i\le \ell$
are linearly independent over a number field $\K$ which satisfies the condition
$\K(\zeta_{\varphi(p_1\cdots p_{\ell})}) \cap \Q(\zeta_{p_1\cdots p_{\ell} }) = \Q$.
\end{cor}

\begin{rmk} 
Consider the sets 
$$
A = \Big\{ ( \frac{p-1}{2},  ~p ) ~\Big{|}~  \text{ both } \frac{p-1}{2} \text{ and } p \text{ are primes} \Big\}
\phantom{m}\text{and}\phantom{m}
B= \Big\{ p ~\Big{|}~ (\frac{p-1}{2},  ~p) \in A \Big\}.
$$
Any prime pair in the set $A$ is called a Sophie-Germain prime pair.
Dickson's conjecture (see preliminaries for precise statement) 
implies the existence of infinitely 
many Sophie-Germain prime 
pairs (see \cite{DEG}).  
Let
$$
C = \Big\{ p_i  ~~\Big{|}~  ~i \ge 1, ~p_i \in B, ~  p_{i+1} >  2p_i + 1  \Big\}
$$
Since by Dickson's conjecture $A$ is an infinite set, so is $B$ and hence $C$ is i
an infinite set.  Choose Erd\"{o}sian functions $f_p$ for $p \in C$.
Then the numbers in the set $\{ L(1, f_{p})  ~|~  p \in C \}$ are linearly
independent over any Galois number field $\K$ whose discriminant $d_{\K}$
is co-prime to $\{ p\varphi(p) ~|~ p \in C\}$. Note that $\K(\zeta_{\varphi(q)})\cap \Q(\zeta_{q}) = \Q$
if and only if $\K \cap \Q(\zeta_{q }) = \Q$
and $\K(\zeta_{\varphi(q)}) \cap \K(\zeta_{q}) = \K$ 
(see \propref{lf} in preliminaries), where $q>2$ is a natural number. Further, the property 
$\K(\zeta_{\varphi(q)}) \cap \K(\zeta_{q}) = \K$ is not necessarily true for $(q, \varphi(q))=1$. 
But when $\K$ is Galois number field whose discriminant $d_{\K}$
is co-prime to $q\varphi(q)$, where $(q, \varphi(q))=1$, then 
$\K(\zeta_{\varphi(q)}) \cap \K(\zeta_{q}) = \K$ (see
Theorem 1.8 in \cite{GMR}).
\end{rmk}

\smallskip 
 
\section{\bf{Preliminaries}}

\smallskip 

In this section, we state the results which will play an important role
in proving our main theorems. We start with the following non-vanishing 
result of Baker, Birch and Wirsing \cite{BBW} (see also chapter 23 of \cite{MR}).

\begin{thm}{\rm (Baker, Birch and Wirsing).}\label{non-vanish}
Let $f$ be a non-zero algebraic valued periodic function with period $q$.
Also let $f(n)=0$ whenever $1 < (n,q) < q$ and
the $q$-th cyclotomic polynomial $\Phi_q(X)$ 
be irreducible over $\Q(f(1),\cdots,f(q))$, then
$$
\displaystyle \sum_{n=1}^{\infty} \frac{f(n)}{n} ~~\neq 0.
$$
\end{thm}

Chowla \cite{chowla} proved that if $p$ is an odd prime, then the
numbers 
$$
\cot(2\pi a/p),  \phantom{m} 1 \le a \le (p-1)/2
$$ 
are linearly independent 
over the field of rational numbers. This result was 
reproved by various authors (see, for instance, \cite{ HH, JH}).
In 1981,  Okada \cite{TO} (see also Wang \cite{KW}) extended Chowla's theorem to 
natural number $q>2$ which are not necessarily primes. In the same theorem
he also considered derivatives of higher orders of $\cot x$. 
Both Okada and Wang made use of the fact that $L(k,\chi)~\ne~0$ 
though their proofs were different.  More precisely, Okada \cite{TO} 
proved the following theorem.
\begin{thm}\label{okada1} 
 Let $k$ and $q$ be positive integers with $k > 0$ and $q > 2$. Let $T$ be a
set of $\varphi(q)/2$ representatives mod $q$ such that the union 
$T \cup ( -T )$ is a complete set
of co-prime residues modulo $q$. Then the set of real numbers
$$
\frac{d^{k-1}}{dz^{k-1}} ( \cot \pi z)|_{(z = a/q)},~ ~ ~  a \in T
$$
is linearly independent over $\Q$.
\end{thm}
Five years later, Girstmair \cite{ KG} gave a much 
simpler proof of this result of Okada using Galois Theory in the case 
when order of the derivative of $\cot x$ is at least $1$.  
In 2009, Murty and Saradha \cite{MS1} extended the work of Okada
to show the following theorem.

\begin{thm}\label{okada}
 Let $k$ and $q$ be positive integers with $k > 0$ and $q > 2$. Let $T$ be a
set of $\varphi(q)/2$ representatives mod $q$ such that the union $T \cup ( -T )$ is a complete set
of co-prime residues modulo $q$. Let $\K$ be an algebraic number field over which the $q$th
cyclotomic polynomial is irreducible. Then the set of real numbers
$$
\frac{d^{k-1}}{dz^{k-1}} ( \cot \pi z)|_{(z = a/q)},~ ~ ~  a \in T
$$
is linearly independent over $\K$.
\end{thm}
See the recent work of Hamahata \cite{YH} for a multi-dimensional generalization
of \thmref{okada}. We deduce another generalization of \thmref{okada} required
for our work.  For this, we will work with linearly disjoint number fields.

\begin{defn} 
Let $\K$ and $\F$ be algebraic extensions of a field $\L$. The fields $\K, \F$ 
are said to be linearly disjoint over $\L$ if every finite subset 
of $\K$ that is $\L$ linearly independent is also $\F$  linearly independent. 
\end{defn}
The following theorem is an equivalent criterion for linearly disjoint fields.

\begin{thm}\label{kf} \cite[Ch 5, Thm 5.5]{MC} 
Let $\K$ and $\F$ be algebraic extensions of a field $\L$. Also let
at least one of $\K, \F$ is separable and one (possibly the same) is normal. 
Then $\K$ and $\F$ are linearly disjoint over $\L$ if and only if $\K \cap \F = \L$.
\end{thm}

We shall also use the following equivalent criterion for linearly disjoint fields.  
\begin{prop}\label{lf} \cite[Ch 5, Prop 5.2]{MC} 
Let $\L \subset \K$ and $\L \subset \E \subset \F$ be algebraic extensions of a field $\L$. 
Then $\K$ and $\F$ are linearly disjoint over $\L$ if and only if $\K$ 
and $\E$ are linearly disjoint over $\L$ and $\K\E$ and $\F$ are linearly disjoint over $\E$.
\end{prop}

In order to prove \thmref{th3},  we use Baker's seminal work on linear forms in logarithms 
of algebraic numbers.

\begin{thm} \label{AB}\cite[Thm 2.1]{AB}
If $\alpha_1,\alpha_2, \cdots, \alpha_n $ are non-zero algebraic numbers such that 
$\log\alpha_1,\log\alpha_2,\\  \cdots, \log\alpha_n $ are linearly independent over the rationals, then $1,\log\alpha_1,\log\alpha_2, \cdots, \log\alpha_n $ are linearly independent over 
the field of algebraic numbers.
\end{thm}

One of the application of the above theorem is the following result which we
will require for our work.

\begin{lem}\label{baker}
\cite[p. 154, Lem 25.4]{MR}  
Let $\alpha_1,\alpha_2, \cdots, \alpha_n $ be positive algebraic numbers. 
If $c_0,c_1, \cdots, c_n$ are algebraic numbers with $c_0 \ne 0$, then 
$$
c_0 \pi + \sum_{j=1}^{n} c_j \log \alpha_j
$$
is a transcendental number.
\end{lem}

For an integer $q > 4$, Ramachandra \cite{KR} discovered a set of multiplicatively independent units 
in the cyclotomic field $\Q(\zeta_q)$, where $\zeta_q$ is a primitive $q$th root of unity. 
For $1< a < q/2$ and $(a,q)=1$, define
$$
\xi_a= \zeta_q^{d_a} \eta_a ~ ~ ~ \in \Q(\zeta_q + \zeta_q^{-1}),
$$
where
$$
d_a = \frac{1}{2} (1-a) \sum_{d \mid q, ~d\ne q \atop{  (d, \frac{q}{d})=1 }} d, \hspace{0.7cm}
\eta_a = \prod_{d \mid q, ~d\ne q \atop{ (d, \frac{q}{d} )=1 } } \frac{1- \zeta_q^{ad}}{ 1 - \zeta_q^d} .
$$
It is easy to see that $\xi_a$ is a unit in $\Q(\zeta_q + \zeta_q^{-1})$ for $1< a < q/2$ and $(a,q)=1$.
Ramachandra proved the following important theorem about these units. 
\begin{thm}\label{units}\cite{KR, LW}
The set of real units
 $\{\xi_a ~ |~ ~  1< a < q/2,(a,q)=1\}$
is multiplicatively independent.
\end{thm}
These units are now known as Ramachandra units. Using these units, one can
express $L(1,\chi)$ when $\chi$ is an even non-trivial character with period $q$
as follows.
\begin{lem}\label{achi}\cite[p. 149]{ LW}
For a natural number $q >4$, let $\chi$ be an even non-trivial character with period $q$.
Then we have
$$
L(1,\chi)  = \delta_{\chi}
\sum_{ 1 < a < q/2 \atop (a,q)=1} \overline{\chi}(a) \log \xi_a,
$$
where $\delta_\chi$ is a non-zero algebraic number. Further, $L(1, \chi)$ can also be written
as algebraic linear combination of logarithms of positive algebraic numbers.
\end{lem}

Using \lemref{achi}, Ram Murty and Kumar Murty \cite{MM1} proved the following theorem.

\begin{thm}\label{f(0)}\cite[Thm 8]{MM1}
Let $q> 2$ be a natural number and $f$ be a non-zero Dirichlet type function with period $q$.
Write $f = f_e + f_o$, where $f_e$ is an even function and $f_o$ is an odd function. 
Let $\K$ be the field generated by the values of $f_o$ over $\Q$. 
If $\K \cap \Q(\zeta_q) = \Q$, then $L(1,f) \ne 0$.
\end{thm}

We end this section by recalling a group theoretic pre-requisite \cite{MM1} as well as 
a conjecture of Dickson \cite{DIC}.
 
\begin{lem}\label{group}
Let $G$ be a finite group. Suppose that for all $g \in G, ~g \ne 1$, we have
$$
\sum_{\chi \ne 1, \atop { \chi} \text{ irreducible}} \mu_{\chi}\chi(g) = 0, \phantom{m}  \mu_{\chi} \in \C,
$$
where the summation varies over all non-trivial irreducible characters 
of $G$. Then $\mu_{\chi} = 0$ for all $\chi \ne 1$.
\end{lem}

\begin{conj}[Dickson's conjecture]\label{dickson}
Let $s$ be a positive integer and  $F_1, F_2,\ldots, F_s$ be $s$ linear polynomials 
with integral coefficients and positive leading coefficient such that their product has no 
fixed prime divisor \footnote{We say that the prime number $p$ is a 
\emph{fixed prime divisor} of a polynomial $G$ if we have: $\forall t \in \mathbb{Z} : p | G(t)$.}. 
Then there exist infinitely many 
positive integers $t$ such that $F_1(t), F_2(t),\ldots, F_s(t)$ are all primes.
\end{conj}

\smallskip

\section{\bf{Proofs of the Main Theorems}}

For $1 \leq j \leq \ell $ and $q_j >2$, consider the sets
$$
S_j = \{ 1 < a_j < q_j /2 ~ | ~ (a_j,q_j) = 1 \}
\phantom{m}\text{and}\phantom{m}
T_j = \{ 1\leq a_j < q_j /2 ~ | ~ (a_j,q_j) = 1 \}
$$
Throughout this section, we shall be using these 
notations.

\subsection{Proof of \thmref{th3}}
We first show that the set of Ramachandra units
$$
\bigcup_{1 \leq j \leq \ell} \{ \xi_{a_j} ~ | ~  a_{j} \in S_j\}
$$
is multiplicatively independent. For $\ell = 1$, it follows from the work
of Ramachandra (see \thmref{units}). Now suppose that 
$$
\bigcup_{1 \leq j < \ell} \{ \xi_{a_j} ~ | ~  a_{j} \in S_j \}
$$
is multiplicatively independent. If there exist $\alpha_{a_{j}} \in \Z$ for
$a_j \in S_j, ~1 \leq j \leq \ell$ such that 
$$
\prod_{1 \leq j \leq \ell} \prod_{a_j \in S_j} \xi_{a_j}^{\alpha_{a_{j}}} = 1,
$$
then
\begin{equation}\label{t-2}
\prod_{1 \leq j < \ell} \prod_{a_j \in S_j} \xi_{a_j}^{\alpha_{a_{j}}} 
=   \prod_{a_{\ell} \in S_{\ell}} \xi_{a_{\ell}}^{ - \alpha_{a_{\ell}}}.
\end{equation}
Note that
\begin{equation}\label{norm}
\prod_{1 \leq j < \ell} \prod_{a_j \in S_j} \xi_{a_j}^{\alpha_{a_{j}}} 
=   \prod_{a_{\ell} \in S_{\ell}} \xi_{a_{\ell}}^{ - \alpha_{a_{\ell}}}
~\in~ 
\Q(\zeta_{q_1 \cdots q_{\ell - 1}}) \cap \Q(\zeta_{q_{\ell}})
~=~ \Q.
\end{equation}
Let us call this rational number $\beta$. If
$\F =\Q(\zeta_{q_1 \cdots q_{\ell}})$ and $N_{\F/\Q}(\alpha)$
denotes the norm of $\alpha$ of $\F$ over $\Q$, then
taking $N_{\F/\Q}$ of the quantities on both sides 
of \eqref{t-2}, we get that $\beta^{ \varphi(q_1\cdots q_{\ell})} =1$
as $\varphi(q_1\cdots q_{\ell})$ is even. This implies that
$\beta = \pm 1$.  Thus
\begin{equation*}
\prod_{1 \leq j < \ell} \prod_{a_j \in S_j} \xi_{a_j}^{ 2\alpha_{a_{j}}} 
~=~   
\prod_{ a_{\ell} \in S_{\ell}} \xi_{a_{\ell}}^{ - 2 \alpha_{a_{\ell}}} 
~=~ 1 .
\end{equation*}
Applying induction hypothesis, we obtain
$\alpha_{a_{j}} = 0$ for  $a_j \in S_j, 1 \leq j \leq \ell$.
This implies that the set of real numbers 
$\bigcup_{ 1 \leq j \leq \ell} \{ \xi_{a_j} ~ | ~ a_{j} \in S_j\}$
 is multiplicatively independent.
 
We now apply the above observation to complete the proof of \thmref{th3}.
Let 
$$
C_j = \{\chi_j  \mod{q_j} ~ | ~ \chi_j(-1) = 1, \chi_j \ne 1 \} .
$$
be the set of non-trivial even characters with periods $q_j$ for $1 \leq j \leq \ell$.
Suppose that there exist algebraic numbers $\alpha_{\chi_{j}}$
 for $\chi_j \in C_j,  1 \leq j \leq \ell$ such that
\begin{equation}\label{eq-rama}
\sum_{ 1 \leq j \leq \ell} \sum_{\chi_j \in C_j} \alpha_{\chi_{j}} L(1,\chi_j) = 0.
\end{equation}
Substituting (see \lemref{achi}) 
$$
L(1,\chi_j)  = \delta_{\chi_{j}} \sum_{a_j \in S_j} \overline{\chi_j}(a_j) \log \xi_{a_{j}}
$$
for $\chi_j \in C_j, 1 \leq j \leq \ell$ in \eqref{eq-rama}, we obtain
\begin{align*}
\begin{split}
\sum_{ 1 \leq j \leq \ell} \sum_{ a_j \in S_j} \left(\sum_{\chi_j \in C_j} \alpha_{\chi_{j}} \delta_{\chi_{j}}  \overline{\chi_j}(a_j)\right) \log \xi_{a_{j}} &= 0.
\end{split}
\end{align*}
Applying Baker's theorem (\thmref{AB}) and our observation about linear independence
of Ramachandra units for $q_1, \cdots ,  q_{\ell}$, we get
$$
\sum_{\chi_j \in C_j} \alpha_{\chi_{j}} \delta_{\chi_{j}}  \overline{\chi_j}(a_j) = 0, 
$$
for $a_{j} \in S_j, 1 \leq j \leq \ell$. Since  $\delta_{\chi_{j}} \ne 0$ (see \lemref{achi}), 
the even characters of $(\Z/q_j\Z)^{\times}$ can be viewed as
characters of the quotient group $(\Z/q_j\Z)^{\times} / \{ \pm1\}$.
As these characters are of dimension one and hence irreducible,
applying \lemref{group}, we have 
$\alpha_{\chi_{j}} = 0$
for $\chi_j \in C_j,  1 \leq j \leq \ell$.
This completes the proof of \thmref{th3}.

\subsection{Proof of \thmref{th4}}
We first show that the set of real numbers
$$
\bigcup_{ 1 \leq j \leq \ell} \{ \cot (\frac{\pi a_j}{q_j}) ~ | ~ a_j \in T_j\}
$$
is linearly independent over $\Q$.
For $\ell = 1$, it follows from the work of Okada (see \thmref{okada1}). 
Suppose that the set of real numbers 
$$
\bigcup_{ 1 \leq j < \ell} \{ \cot (\frac{\pi a_j}{q_j}) ~ | ~ a_j \in T_j\}
$$
is linearly independent over $\Q$.  If there exist rational numbers 
$\alpha_{a_{j}}$ for $a_j \in T_j, 1 \leq j \leq \ell$ such that 
\begin{align*}
\begin{split}
\sum_{ 1 \leq j \leq \ell} \sum_{a_j \in T_j} \alpha_{a_{j}} \cot (\frac{\pi a_j}{q_j}) &= 0,
\end{split}
\end{align*}
then
\begin{align}
\begin{split}\label{t-1}
\sum_{ 1 \leq j < \ell} \sum_{a_j \in T_j} \alpha_{a_{j}} \cot (\frac{\pi a_j}{q_j}) &= - \sum_{a_{\ell} \in T_{\ell}} \alpha_{a_{\ell}} \cot (\frac{\pi a_{\ell}}{q_{\ell}}).
\end{split}
\end{align}
Since 
$$
- i \cot \frac{\pi a_j}{q_j} ~=~ \frac{\zeta^{a_j}_{q_j} + 1}{\zeta^{a_j}_{q_j} - 1} \in \Q(\zeta_{q_j}),
$$  
where $ i = \sqrt{-1}$, it follows that
$$
i \sum_{ 1 \leq j < \ell} \sum_{a_j \in T_j} \alpha_{a_{j}} \cot (\frac{\pi a_j}{q_j}) 
~=~ 
- i \sum_{a_{\ell} \in T_{\ell}} \alpha_{a_{\ell}} \cot (\frac{\pi a_{\ell}}{q_{\ell}})
~\in~ 
\Q(\zeta_{q_1 \cdots q_{\ell - 1}}) \cap \Q(\zeta_{q_{\ell}})
~=~ \Q.
$$
Since a purely imaginary number is a rational number if and only if it is $0$, we have
\begin{align*}
\begin{split}
\sum_{ 1 \leq j < \ell} \sum_{a_j \in T_j} \alpha_{a_{j}}  \cot (\frac{\pi a_j}{q_j}) = - \sum_{a_{\ell} \in T_{\ell}} \alpha_{a_{\ell}} \cot (\frac{\pi a_{\ell}}{q_{\ell}}) &= 0.
\end{split}
\end{align*}
Applying induction hypothesis, we get that $\alpha_{a_{j}} = 0$ for all $a_j \in T_j,  1 \leq j \leq \ell$. 
Hence the set of real numbers 
$$
 \bigcup_{ 1 \leq j \leq \ell} \{ \cot (\frac{\pi a_j}{q_j}) ~ | ~ a_j \in T_j\}
 $$
 is linearly independent over $\Q$.

We now apply the above observation to complete the proof of \thmref{th4}.
Let 
$$
D_j = \{\chi_j \mod{q_j} ~ | ~ \chi_j(-1) = -1\}
$$
be the set of odd characters with periods $q_j$ for $1 \leq j \leq \ell$. Let $\K$ be as in \thmref{th4}.
Suppose that there exist $\alpha_{\chi_{j}} \in \K(\zeta_{\varphi(q_1) \cdots \varphi(q_{\ell})})$ 
for $\chi_j \in D_j, 1 \leq j \leq \ell$ such that 
\begin{equation}\label{der}
\sum_{ 1 \leq j \leq \ell} \sum_{\chi_j \in D_j} \alpha_{\chi_{j}} L(1,\chi_j) = 0.
\end{equation}
Substituting (see \cite{{NR}, {TO}})
\begin{equation}\label{odd-chi}
L(1,\chi_j)  = \frac{\pi}{q_j} \sum_{a_j \in T_j} \chi_j(a_j) \cot (\frac{\pi a_j}{q_j}), 
\end{equation}
for $\chi_j \in D_j,1 \leq j \leq \ell$ in \eqref{der}, we obtain
\begin{align}
\begin{split}\label{lkchi}
\sum_{ 1 \leq j \leq \ell} \sum_{a_j \in T_j} \left(\sum_{\chi_j \in D_j} 
\frac{\alpha_{\chi_{j}}}{q_j} \chi_j(a_j)\right) \cot (\frac{\pi a_j}{q_j}) &= 0.
\end{split}
\end{align}
By given hypothesis and \thmref{kf}, the number fields 
$\K(\zeta_{\varphi(q_1) \cdots \varphi(q_{\ell})})$ and $\Q(\zeta_{q_1 \cdots q_{\ell}})$ 
are linearly disjoint over $\Q$. Therefore $\Q$-linearly independent elements  
$i\cot (\frac{\pi a_j}{q_j})$ in \eqref{lkchi} which
belong to $\Q(\zeta_{q_1 \cdots q_{\ell}})$ are also linearly independent over
$\K(\zeta_{\varphi(q_1) \cdots \varphi(q_{\ell})})$.
Since the coefficients of $\cot (\frac{\pi a_j}{q_j})$ in~\eqref{lkchi}
belong to $\K(\zeta_{\varphi(q_1) \cdots \varphi(q_{\ell})})$, we have
$$
\sum_{\chi_j \in D_j} \frac{\alpha_{\chi_{j}}}{q_j} \chi_j(a_j) = 0
$$
for $a_j \in T_j ,1 \leq j \leq \ell$. Since all the characters in the set $D_j, 1 \leq j \leq \ell$ are 
of same parity, it follows that 
$$
\sum_{\chi_j \in D_j} \frac{\alpha_{\chi_{j}}}{q_j} \chi_j(a_j) = 0
$$ 
for $a_j \in (\Z/q_j\Z)^{\times}  , 1 \leq j \leq \ell$. It then follows
from linear independence of characters that
$\alpha_{\chi_{j}} = 0$ for $\chi_j \in D_j, 1 \leq j \leq \ell$.
This completes the proof of \thmref{th4}.

\subsection{Proofs of \thmref{th1} and \thmref{th2}}

Note that \thmref{th1} follows by considering $\K =\Q$ in 
\thmref{th2}. Hence it is sufficient to prove \thmref{th2}.
It follows from \lemref{achi} that for an even non-trivial Dirichlet character $\chi$, 
the number $L(1, \chi)$ is a linear form in logarithms of positive real algebraic numbers.
We know from \eqref{odd-chi} that for an odd character $\chi $,  the number $L(1, \chi)$ 
is an algebraic multiple of $\pi$. Then \lemref{baker} implies that the space 
generated by $L(1, \chi)$ for non-trivial even $\chi$ do not intersect with 
the space generated $L(1, \chi)$ for odd $\chi$.  \thmref{th2} now follows
by applying \thmref{th3} and  \thmref{th4}.

\subsection{Proof of \corref{cor1}}
Let us denote by $Y_{j,e} = \{ f \in Y_{q_j}(\overline{\Q}) ~|~   f(-a) = f(a) \text{ for } 1 \le a <q \}$.
Suppose that there exist  $\alpha_{f_{j}} \in \overline{\Q}$ for 
$f_j \in Y_{j, e} , 1 \leq j \leq \ell$ such that
$$
\sum_{1 \leq j \leq \ell} \sum_{f_j \in Y_{j,e} } \alpha_{f_{j}} L(1,f_j) ~=~ 0.
$$
Then 
\begin{equation}\label{e-11}
\sum_{ 1 \leq j \leq \ell} L(1, F_j) ~=~ 0,
\end{equation}
where $F_j= \sum_{f_j \in Y_{j, e} } \alpha_{f_j} f_j$. For 
$1 \leq j \leq \ell$,  $f_j \in Y_{j,e}$ and hence $F_j$'s are even
Dirichlet type functions with periods $q_j$.
Therefore we can write $F_j$ as a linear combination of $\chi_j$,
where $\chi_j$ belong to the set $C_j = \{\chi  \mod{q_j} ~ | ~ \chi(-1) = 1, \chi \ne 1\} $.
This implies that
$$
L(1, F_j) = \sum_{\chi_j \in C_j} \beta_{\chi_{j}} L(1,\chi_j),
$$
where $\beta_{\chi_{j}}$ are algebraic numbers.
Substituting this
expression in \eqref{e-11}, we get
\begin{align*}
	\begin{split}
		\sum_{ 1 \leq j \leq \ell} \sum_{\chi_j \in C_j}
		\beta_{\chi_{j}} L(1,\chi_j) &= 0.
	\end{split}
\end{align*}
Applying \thmref{th3}, we obtain $\beta_{\chi_{j}} = 0$ for 
$\chi_j \in C_j,   ~1 \leq j \leq \ell$. Thus $L(1, F_j) = 0$ for $1 \leq j \leq \ell$.
Using \thmref{f(0)} (see also \cite[Th. 6]{MM1}), we then have $F_j= 0$ 
for $1 \leq j \leq \ell$. Since by hypothesis, the elements of $Y_{j, e}$ are $\overline{\Q}$
linearly independent, we have $\alpha_{f_{j}} = 0$
for $ f_j \in Y_{j,e} , ~1 \leq j \leq \ell$.
This completes the proof of \corref{cor1}.

\subsection{Proof of \corref{cor3} }

Let $\{ f_p ~|~ p \text{ odd prime}\}$ be as in \corref{cor3}. 
Note that we can write $f_p$ as a sum of an even function
and an odd function, i.e., $f_p = f_{p,e} + f_{p,o}$, 
where 
$$
f_{p,e}(a) ~=~\frac{f_p(a) + f_p(-a)}{2}
\phantom{m}\text{and}\phantom{m} 
f_{p,o}(a) ~=~ \frac{f_p(a) - f_p(-a)}{2}
$$
for $1\le a \le q$. If the corollary is not true, then there exist a finite subset 
$\mathcal{P}$ of prime numbers and algebraic numbers $\alpha_{p}$ (not all zero) 
for $p \in \mathcal{P}$ such that 
\begin{equation}\label{sup}
\sum_{ p \in \mathcal{P}} \alpha_{p} L(1, f_p) ~=~ 0.
\end{equation}
This implies that
$$
\sum_{ p \in \mathcal{P}} \alpha_{p} L(1,f_{p,e}) 
~+~ 
\sum_{ p \in \mathcal{P}} \alpha_{p} L(1,f_{p,o}) 
~=~ 0.
$$
Since each $L(1,f_{p,e})$ for $p \in \mathcal{P}$ can be written as
algebraic linear combination of $L(1, \chi)$'s for non-trivial
even Dirichlet characters $\chi$ with period $p$, it follows from
\lemref{achi} that the summation 
 $\sum_{ p \in \mathcal{P}} \alpha_{p} L(1,f_{p,e})$
 is an algebraic linear combination of logarithms of positive algebraic 
numbers. Similarly each $L(1,f_{p,o})$ for $p \in \mathcal{P}$ can be written as
algebraic linear combination of $L(1, \chi)$'s for odd
Dirichlet characters $\chi$ with period $p$, we see that
$\sum_{ p \in \mathcal{P}} \alpha_{p} L(1,f_{p,o}) $
is an algebraic multiple of $\pi$
by identity \eqref{odd-chi}. Now by applying \lemref{baker}, we have
$$
\sum_{ p \in \mathcal{P}}\alpha_{p} L(1,f_{p,e}) 
~=~ 
- \sum_{ p \in \mathcal{P}} \alpha_{p} L(1,f_{p,o})
~=~ 0.
$$
Since $f_{p,e}$ are non-zero even Dirichlet type functions with 
distinct prime periods $p \in \mathcal{P}$, we have 
$L(1, f_{p,e})$ are non-zero for $p \in \mathcal{P}$.
Now applying \corref{cor1}, we have
$\alpha_p = 0$ for $p \in \mathcal{P}$, a contradiction
to \eqref{sup}. This completes the proof of \corref{cor3}.

\subsection{ Proof of Corollary \ref{cor2}}
Let $\K$ be as in \corref{cor2} and for $1 \le j \le \ell$, $Y_j(\K)$ denotes $Y_{q_j}(\K)$ 
for the sake of brevity.  As in \corref{cor3}, let us write
$f_j = f_{j,e} + f_{j, o}$, where
$$
f_{j,e}(a) ~=~\frac{f_j(a) + f_j(-a)}{2}
\phantom{m}\text{and}\phantom{m} 
f_{j,o}(a) ~=~ \frac{f_j(a) - f_j(-a)}{2}
$$
for $1\le a \le q$ and $1 \le j \le \ell$. 
Suppose that there exist $\alpha_{f_{j}} \in \K$ for 
$f_j \in Y_j(\K) , 1 \leq j \leq \ell$ such that 
$$
\sum_{ 1 \leq j \leq \ell} 
\sum_{f_j \in Y_j(\K) } \alpha_{f_{j}} L(1, f_j) ~=~ 0.
$$
This implies that
\begin{equation}\label{eq-eo}
\sum_{ 1 \leq j \leq \ell} 
\sum_{f_j \in Y_j(\K)} \alpha_{f_{j}} L(1,f_{j,e}) 
~+~ 
\sum_{ 1 \leq j \leq \ell} \sum_{f_j \in Y_j(\K)}
\alpha_{f_{j}} L(1,f_{j,o}) 
~=~ 0.
\end{equation}
Proceeding as in \corref{cor3}, we note that the first term in \eqref{eq-eo} 
is an algebraic linear combination of logarithms of positive algebraic numbers
by \lemref{achi}  and the second term of~\eqref{eq-eo} 
is an algebraic multiple of $\pi$ by identity \eqref{odd-chi}. Applying \lemref{baker}, we have
\begin{equation*}
\sum_{ 1 \leq j \leq \ell} \sum_{f_j \in Y_j(\K)} \alpha_{f_{j}} L(1,f_{j,e}) 
~=~  0 
\phantom{m}\text{and}\phantom{m}
 \sum_{ 1 \leq j \leq \ell} \sum_{f_j \in Y_j(\K)} \alpha_{f_{j}} L(1,f_{j,o}) 
 ~=~  0.
\end{equation*}
This implies that 
\begin{equation}\label{e-15}
\sum_{1 \leq j \leq \ell} L(1, F_{j, e}) ~=~ 0
\phantom{m}\text{and}\phantom{m}
\sum_{1 \leq j \leq \ell} L(1, F_{j, o}) ~=~ 0,
\end{equation}
where
\begin{equation}
F_{j, e} 
~=~ 
\sum_{f_j \in Y_j(\K) }\alpha_{f_{j}} f_{j,e}
\phantom{m}\text{and}\phantom{m}
F_{j, o} 
~=~ 
\sum_{f_j \in Y_j(\K) }\alpha_{f_{j}} f_{j,o}.
\end{equation}
Since $F_{j, e}$'s are even Dirichlet type functions with
distinct periods $q_j$ for $1 \le j \le \ell$, applying \corref{cor1},
we have 
\begin{equation}\label{eq-ez}
F_{j, e} ~=~0
\phantom{m}\text{ for } 1 \le j \le \ell.
\end{equation}
Note that $F_{j, o}$'s are odd Dirichlet type functions with periods $q_j$ for
$1 \le j \le \ell$ with values in $\K$.
Let $V_j$ be the $\K(\zeta_{\varphi(q_j)} )$ vector space of 
functions from $(\Z/q_j\Z)^{\times}$ to $\K(\zeta_{\varphi(q_j)})$.  
Dirichlet characters with periods $q_j$ are contained in $V_j$ and they form 
a basis of $V_j$ over $\K(\zeta_{\varphi(q_j)}) $. Let $D_j$ be the set of all odd
Dirichlet characters with periods $q_j$. Since $F_{j, o} $ can be 
written as
$$
F_{j, o} = \sum_{\chi_j \in D_j} \beta_{\chi_j} \chi_j
$$ 
where $\beta_{\chi_j}  \in \K(\zeta_{\varphi(q_j)} )$ for $\chi_j \in D_j,   1 \leq j \leq \ell$,
we have
$L(1, F_{j, o}) = \sum_{\chi_j \in D_j} \beta_{\chi_{j}} L(1,\chi_j)$. 
Substituting this expression in \eqref{e-15}, we have
\begin{align*}
	\begin{split}
		\sum_{ 1 \leq j \leq \ell} \sum_{\chi_j \in D_j}
		\beta_{\chi_{j}} L(1,\chi_j) &= 0.
	\end{split}
\end{align*}
Since by hypothesis,
$\K(\zeta_{\varphi(q_1) \cdots \varphi(q_{\ell})}) \cap \Q(\zeta_{q_1 \cdots q_{\ell}}) = \Q$, 
applying \thmref{th4}, we get $\beta_{\chi_{j}} = 0$
for $\chi_j \in D_j,  ~1 \leq j \leq \ell$. This implies that
$$
L(1, F_{j, o}) = 0
$$ 
for $1 \leq j \leq \ell$. \thmref{f(0)} then implies that
\begin{equation}\label{eq-ez}
F_{j, o} ~=~0
\phantom{m}\text{ for } 1 \le j \le \ell.
\end{equation}
Then for $1 \le j \le \ell$, we have
$$
\sum_{f_j \in Y_j(\K) }\alpha_{f_{j}} f_{j}
~=~ 
\sum_{f_j \in Y_j(\K) }\alpha_{f_{j}} (f_{j,e} + f_{j,o})
~=~
F_{j,e} + F_{j, o} 
~=~0.
$$
Since by hypothesis, elements of $Y_j(\K) $ are
$\K$ linearly independent, we have 
$\alpha_{f_{j}} =0$ for any $f_j~\in Y_j(\K) , 1 \leq j \leq \ell$.
This completes the proof of \corref{cor2}.

\bigskip
\noindent
{\bf Acknowledgments.} 
The authors thank Purusottam Rath for going through an earlier
version of the article and for many valuable comments.
The first author acknowledges MTR/2018/000201
and SPARC project~445 for partial financial support.
Both the authors would like to thank DAE number 
theory plan project.


\medskip

\begin{thebibliography}{100}
\bibitem{AB}
A. Baker,
Transcendental number theory, 
Cambridge University Press, {\em Cambridge}, 1975.


\bibitem{BBW} 
A. Baker,  B.J. Birch and  E. A Wirsing,
{\em On a problem of Chowla},
J. Number Theory, {\bf 5} (1973), 224--236. 


\bibitem{chowla}
 S. Chowla,
 {\em The nonexistence of nontrivial linear relations between the roots of 
 a certain irreducible equation}, 
 J. Number Theory, {\bf 2} (1970), no. 2, 120--123.
 

\bibitem{MC}
P.  M. Cohn, 
Algebra, Second edition, vol. 3, {\em John Wiley \& Sons}, 1991.


\bibitem{DEG} 
J. M. Deshouillers, P. Eyyunni and S. Gun,  
{\em On the local structure of the set of values of Euler's $\varphi$ function},
Acta Arith., {\bf 199} (2021), no. 1, 103--109. 
 
 
\bibitem{DIC} 
 L. E. Dickson, 
 {\em A new extension of Dirichlet's theorem on prime numbers}, 
 Messenger of Math., {\bf 33} (1904), 155--161.

 
\bibitem{KG}  
K. Girstmair, 
{\em Letter to the editor}, 
J. Number Theory, {\bf 23} (1986), p. 405.


\bibitem{GMR}
S. Gun, M. Ram Murty and P. Rath,
{\em Linear independence of Hurwitz zeta values and a theorem of 
Baker-Birch-Wirsing over number fields},
 Acta Arith. {\bf 155} (2012), no. 3, 297--309.


\bibitem{YH}
Y. Hamahata, 
{\em Okada's theorem and multiple Dirichlet series}, 
Kyushu J. Math, {\bf 74} (2020), 429--439.


\bibitem{HH}
H. Hasse, 
{\em On a question of S. Chowla},
Acta Arithmetica,  {\bf 18} (1971), 275--280.


\bibitem{JH}
H. Jager and H.W. Lenstra,
{\em Linear independence of cosecant values}, 
Nieuw Arch. Wisk, {\bf 23} (1975), no. 3, 131--144.


\bibitem{NR} 
N. Kandhil and P. Rath, 
{\em Around a question of Baker}, submitted.


\bibitem{MS1}  
M. Ram Murty and N. Saradha,
{\em Special values of the polygamma functions,} 
Int. J. Number Theory, {\bf 5} (2009), no. 2, 257--270. 


\bibitem{MM1}  
M. Ram Murty and V. Kumar Murty,
{\em A problem of Chowla revisited,} 
J. Number Theory, {\bf 131} (2011), no. 9, 1723--1733.


\bibitem{MR} 
 M. Ram Murty and P. Rath,
Transcendental numbers, 
{\em Springer}, New York, 2014.


\bibitem{KR} 
K. Ramachandra, 
{\em On the units of cyclotomic fields,}
Acta Arith, {\bf 12}, (1966/67), 165--173.


\bibitem{TO} 
T. Okada, 
{\em On an extension of a theorem of S. Chowla,}  
 Acta Arithmetica, {\bf 38} (1981), no. 4, page 341--345.


\bibitem{KW}
K. Wang, 
{\em On a theorem of S. Chowla},
J. Number Theory, {\bf 15} (1982), 1--4.


\bibitem{LW}
L. C. Washington, 
Introduction to cyclotomic fields,
Second edition, Graduate Texts in Mathematics, 
{\em Springer-Verlag}, New York, 1997.
 

\end{thebibliography}
\end{document}